\documentclass{conm-p-l}
\usepackage{latexsym}
\usepackage{amsfonts}
\usepackage{graphicx}
\usepackage{amssymb}
\usepackage{amsmath}

\title[Pattern Recognition in Free Groups]{Pattern Recognition Approaches to Solving Combinatorial
Problems in Free Groups}
\author[R. M. Haralick]{Robert M. Haralick}
\address{Department of Computer Science, Graduate Center,
City University of New York, 365 Fifth Avenue, New York, NY 10016}

\author[A. D. Miasnikov]{Alex D. Miasnikov}

\author[A. G. Myasnikov]{Alexei G. Myasnikov}

\newcommand{\TwoMap}[2]{ \left ( \begin{array}{ll} a \rightarrow #1 \\ b
\rightarrow #2 \end{array} \right ) }

\begin{document}

\maketitle

\begin{abstract}
We review some basic methodologies from pattern recognition that
can be applied to helping solve combinatorial problems in free group theory.
We illustrate how this works with recognizing Whitehead minimal words
in free groups of rank 2. The methodologies reviewed include how to form
feature vectors, principal components, distance classifers, linear classifiers,
regression classifiers, Fisher linear discriminant, support vector machines,
quantizing, classification trees, and clustering techniques.
\end{abstract}

\keywords{Free group,  automorphism problem, Whitehead method,
Pattern Recognition, Classification}

\subjclass{Primary 20F28, Secondary 68T10}

\tableofcontents

\section{Introduction}
There are many problems in mathematics that can be solved by a
search over all possibilities. In combinatorial group theory every
classical decision problem has its natural "search" variation. For
example, the search  version of the Word Problem for a group given
by a  presentation  $\langle X \mid R \rangle$ would ask for a
word $u$ from the normal closure of the set of relators $R$ to
produce an expression of $u$ as a product of conjugates of
elements from $R^{\pm 1}$. The search  Conjugacy Problem would
require to produce a conjugating element, and the search
Membership Problem would ask for an expression of a given $u$ from
a subgroup $\langle v_1, \ldots, v_k\rangle$ as a product of the
generators (and their inverses) $v_1^{\pm 1}, \ldots, v_k^{\pm
1}$. In a free group finding a solution of a given consistent
equation or determining whether or not  a given element has
minimal length in its automorphic orbit (Whitehead minimality
problem)  are typical search problems. All these problems are
recursively enumerable, so,  in principal, the  total search would
produce an answer sooner or later. Unfortunately, in practice  the
total search could be extremely inefficient: there are  not any
recursive bounds on time complexity of the search word, conjugacy,
or membership problems, and  there is no algorithm known to find
solutions of equations or determine minimality of a word in a free
group in time better than exponential time in the size of the
problem. In this paper, we offer some insights to speed up the
solving of computational hard problems that arise in certain
problem populations of  free groups. We will show that with some
intelligent reasoning, a significant fraction of problems in the
population can be solved much quicker than exponential time.

The classical way algorithms are formulated in search problems is
to design a way to perform the associated tree search and try to
improve the solution time by fancier data structures, more
efficient code, and by heuristics. Our point of view is to improve
how we solve these problems using the experience of having solved
many such problems from the given population. Thus our problem
domain is given by a population of problems. We will sample many
problem instances from the population and solve each instance. We
will then examine the statistical characteristics in the tree
search of how each problem instance was solved and then
incorporate this statistical derived knowledge from the experience
of solving these problems into a smarter tree search. There are
situations in which we may not have to do any tree searching at
all and our pattern recognition techniques will do all the work.

The use of this kind of experience may not improve the worst case
complexity. But it can dramatically improve the computational
complexity on a large fractions of problems from the population.
For example, we may discover that 99 percent of the problems can
be solved in linear time, 0.9 percent in quadratic time, and the
remainder 0.1 percent in exponential time.

There are two dimensions to using this pattern recognition
technology. The first dimension is the representation of any
problem instance to a fixed dimensional feature vector which
captures the information in the problem instance. The second
dimension is the use of standard tools in pattern recognition that
given the feature vector designates the most probable class or the
most probable next successful step in the tree search.

In this paper we discuss briefly several typical pattern
recognition techniques and demonstrate some of their applications
by the Whitehead's minimality problem. {\em For notations and
relevant results on Whitehead method we refer to the paper}
\cite{G:Miasnikov03GAWhitehead}.

\section{General remarks on Pattern Recognition tasks}
\label{sub-se:1}

One of the main applications of Pattern Recognition ({\bf PR})
techniques is classification of a variety of given objects into
categories. Usually classification algorithms or
\emph{classifiers} try to find a set of measurements (properties,
characteristics) of objects,   called \emph{features}, which gives
a  descriptive representation for the objects.

 Generally, pattern recognition techniques can be divided in two principal
types:
\begin{itemize}
\item \emph{ supervised learning; } \item \emph{ unsupervised
learning {\rm (}clustering\/{\rm )}. }
\end{itemize}
In supervised learning the decision algorithms are ``trained'' on
a prearranged dataset, called \emph{training} dataset in which
each pattern is labelled with its true class label.   If such
information is not available, one can use clustering. In this case
clustering algorithms try to find \emph{clusters} or ``natural
groupings'' of the given objects. In this paper we use supervised
learning pattern recognition algorithms.

Every pattern recognition task of the supervised learning type has
to face all of the following issues:
\begin{enumerate}
 \item {\bf Obtaining the data. } The {\em training datasets}  can be
obtained from the real world or generated by reliable  procedures,
which provide   independent and representative  sampled data.

 \item {\bf Feature extraction. }  The task of feature extraction is
problem specific and requires knowledge of the problem domain. If
such  knowledge  is limited then one may consider as many features
as possible and then try to extract the most ``significant'' ones
using statistical methods.

\item {\bf Model selection. }  The model is the  theoretical basis
for the classifier. A particular choice of the model determines
the basic design of the classifier (though there might be some
variations in the implementation). Model selection is one of the
most active areas of research in pattern recognition. Usually
model selection is closely related to the feature extraction. In
practice, one may try several standard models starting with the
simplest ones or more economic ones.

\item {\bf Evaluation. } Evaluation of the performance of a
particular {\bf PR}  system is an important component of the
pattern recognition task.  It answers the question whether or not
the given system performs as required. To evaluate a system one
can use various accuracy measures, for example, percentage of
correct answers.  To get  reliable estimates other sets of data
that are independent from the training sets must be used. Such
sets are called {\em test datasets}.

Typically we view a {\bf PR} system as consisting of components
1)-4).

\item {\bf Analysis of the system. }
 Careful analysis of performance of a particular classifier may improve
feature extraction and model selection. For example, one can look
for an optimal  set of features or for a more effective model.
Moreover, through analysis of the most significant (insignificant)
features one may gain a new knowledge about the original objects.
\end{enumerate}

\section{Feature Vectors}
\label{se:feature}

 A feature vector is just a vector of properties
about the object of interest, in our case about the mathematical
object of interest. The properties are chosen to be ones thought
to be relevant to the problem. We will illustrate the selection of
features by the Whitehead minimal problem for a free group  $F =
F(X)$ with  basis $X$.

Let $w$ be a reduced word in the alphabet $ \in X^{\pm 1}$. Below
we describe  features of  $w$ which characterize a certain
placement of specific words  from $F(X)$ in $w$.

  Let $K \in \mathbb{N}$ be a natural number,  $v_1, \ldots, v_{K} \in F(X)$
  be words from $F(X)$, and $U_1, \ldots, U_{K+1} \subseteq F(X)$ be subsets
  of $F(X)$. Denote by
\[C(w,U_1v_1\ldots v_KU_{K+1})\]
 the number of  subwords of the type
 $$ u_1 v_1 u_2 \dots v_K u_{K+1},$$
 where $  u_j \in U_j,$ which occur in $w$.   For fixed $K, v_1, \ldots, v_K,
 U_1, \ldots, U_{K+1}$   we obtain a {\em counting function}
 \begin{equation}
 \label{eq:feature}
 w \in F \longrightarrow
C(w,U_1v_1\ldots v_KU_{K+1}) \in \mathbb{N}
\end{equation}
 The normalized value
 $$\frac{1}{|w|}C(w,U_1v_1\ldots v_KU_{K+1})$$
  is called a {\em feature} of $w$ and the function
  $$w \in F \longrightarrow \frac{1}{|w|}C(w,U_1v_1\ldots v_KU_{K+1}) \in
\mathbb{R}$$ is called a {\em feature function} on $F$. Usually we
omit $U_i$ in our notations if $U_i = \emptyset$.  If ${\bar C}=
(C_1(w), \ldots, C_N(w))$ is a sequence of  counting functions
like (\ref{eq:feature}) one  can
 associate with $w$ a vector of real numbers:
$$f_{\bar C}(w) = \frac{1}{|w|}< C_1(w), \ldots, C_N(w)> \ \in \mathbb{R}^N$$
which is called a {\em feature vector}. Every choice of the
sequence ${\bar C}$ gives a vector $f_{\bar C}(w)$ which reflects
the structure of $w$.

For example, if $a \in X^{\pm 1}$ then $C(w,a)$ counts the number
of occurrences of the letter $a$ in $w$. The feature vector (where
for simplicity  we assume that the components are written in some
order which  we do not specify)
 $$f_0(w) = \frac{1}{|w|}<C(w,a) \mid a \in X^{\pm 1}>$$
 shows the frequencies of occurrences of letters from $X^{\pm 1}$
 in $w$. The  feature vector
 $$f_1(w) = \frac{1}{|w|}<C(w,v) \mid |v| = 2>$$
shows the numbers of occurrences of words of length two in $w$
relative to the length of $w$.

To visualize some  structures described by  the counting functions
above we associate with a given word $w \in F(X)$ a  weighted
labelled directed  graph $\Gamma(w)$. Put $V(\Gamma(w)) = X^{\pm
1}$. For given $x, y \in X^{\pm 1}$ and $v \in F(X)$  we connect
the vertex $x$ to the vertex $y$ by an edge with a label $v$ and
weight $C(w,xvy)$. Now, with every edge from $x$ to $y$ with label
$xvy$ one can associate a  counting function $C(w,xvy)$, and vice
versa. It follows that every subgraph $\Gamma$ of $\Gamma(w)$
gives rise to a particular set of counting functions ${\bar
C}_{\Gamma}$ of the type $C(w,xvy)$, and conversely, every  set
${\bar C}$  of counting functions of the type $C(w,xvy)$
determines a subgraph $\Gamma_{\bar C}$ of $\Gamma(w)$. For
instance, the feature mapping $f_1$ corresponds to the subgraph
$\Gamma_1(w)$ of $\Gamma(w)$ which is in a sense a directed
version of the so-called {\em Whitehead graph} of $w$.

  Let $U_n$ be the
set of all words in $F$ that are length $n$. Let $W_n$ be the set
of all words in $F$ that are of length $n$ or less. Other relevant
features can be defined as follows. Each corresponds to various
subgraphs of the graph $\Gamma(w)$:

\[f_2(w) = \frac{1}{|w|}<C(w,x_1 U_1 x_2) \mid x_1,x_2 \in X^{\pm 1}>;\]
\[f_3(w) = \frac{1}{|w|}<C(w,x_1 U_2 x_2) \mid  x_1,x_2 \in X^{\pm 1}>;\]
\[f_4(w) = \frac{1}{|w|}<C(w,x_1 U_3 x_2) \mid  x_1,x_2 \in X^{\pm 1}>;\]
\[f_{5}(w) = \frac{1}{|w|}<C(w,x_1 W_1 x_2) \mid  x_1,x_2 \in X^{\pm 1}>;\]
\[f_{6}(w) = \frac{1}{|w|}<C(w,x_1 W_3 x_2) \mid  x_1,x_2 \in X^{\pm 1}>.\]

\section{Pattern Recognition Tools and Models}

There are a variety of pattern recognition tools that are useful
for determining a way of making a distinction given a set of
feature vectors from one class and a set of feature vectors from
another class.  For each given class of objects, we sample objects
from the class and construct the set of corresponding feature
vectors. The pattern recognition technology provides a way of
determining a best or near best boundary in the feature space that
distinguishes the one class from the other. In this section we
will review some of the basic techniques. The reader interested in
a fuller discussion may consult general references
\cite{duda},\cite{patrick}\cite{fukunaga},\cite{schalkoff}.

\subsection{Principal Components}
There are occasions when the feature vectors coming from a class
either all lie in a small dimensional flat or most of them lie in
a small dimensional flat. Principal component analysis can
determine this.

Let $x_1,\ldots, x_N$ be the set of feature vectors sampled from a
given class. Define their sample mean $\mu$ by
$$
\mu = \frac{1}{N}\sum_{n=1}^N x_n
$$
Define their sample covariance $C$ by
$$
C=\frac{1}{N}\sum_{n=1}^N (x-\mu)(x-\mu)'
$$
Let $t_1,\ldots,t_N$ be the eigenvectors of $C$ with corresponding
eigenvalues $\lambda_1 \ge \lambda_2 \ldots \ge \lambda_N$. Should
the feature vectors indeed lie in a small dimensional flat, then
there will be a $K < N$ such that $\lambda_k = 0, \ k=K+1,\ldots,
N$.

In this case, feature vectors $x$ coming from objects that are in
the class can be recognized by testing whether or not
$$
|| T(x-\mu)|| = 0
$$
where $T$ is a $N\times (N-K)$ matrix whose columns are
eigenvectors $t_{K+1},\ldots,t_N$.

Those in the class will have $||T(x-\mu)|| = 0$. $||T(x-\mu)|| >
0||$ is a sure indication that $x$ comes from an object out of the
class, but there may be some objects out of the class for which
$||T(x-\mu)|| =0$.

In the case of two classes, we form the matrix $T_1$ from the zero
eigenvalue eigenvectors of the covariance matrix from the class
one feature vectors and the matrix $T_2$ from the zero eigenvalue
eigenvectors of the covariance matrix from the class two feature
vectors.

Now if $||T_1 (x-\mu_1)|| = 0$ and $||T_2 (x-\mu_2)|| > 0$, we
assign vector $x$ to class one. If $||T_1 (x-\mu_1) || > 0$ and
$||T_2 (x-\mu_2)|| = 0$, we assign vector $x$ to class two. If $||
T_1 (x-\mu_1)|| = 0 $ and $||T_2 (x-\mu_2)|| = 0$, feature vector
$x$ comes from an object that is both class 1 and class 2. If
$||T_1 (x-\mu_1)|| > 0$ and $||T_2 (x-\mu_2)|| > 0$ feature vector
$x$ comes from an object that is neither class 1 nor class 2.

\subsection{Classifying by Distance}
Let $T_1$ and $T_2$ be defined as before. We form the discriminant
function
$$
f(x) = || T_1(x-\mu_1)|| - || T_2(x-\mu_2)||
$$
which measures the difference between feature vector $x$ and the
flat associated with class 1 and the flat associated with class 2.
The decision rule is to assign vector $x$ to class 1 if $f(x) >
\theta$, otherwise assign to class 2. Here, after the discriminant
function is defined we determine the value of $\theta$ that
minimizes the error.

Classifying by distance can also be done with respect to the class
means. Here the discriminant function is defined by
$$
f(x) = (x-\mu_1)'C_1^{-1}(x-\mu_1) - (x-\mu_2)'C_2^{-1}(x-\mu_2)
$$
As before, the decision rule is to assign vector $x$ to class 1 if
$f(x) > \theta$, otherwise assign to class 2. Here also after the
discriminant function is defined we determine the value of
$\theta$ that minimizes the error.

\subsection{Linear Classifiers}
\label{se:linear}

Classifying may be done by a linear decision rule. Here the
discriminant function is given by
$$
f(x) = v'x
$$
where vector $v$ is the weight vector and is the normal to the
hyperplane separating the feature space into two parts.

If $f(x) < \theta$ the decision rule is to assign the vector $x$
to class 1 otherwise to class 2.There are a variety of ways to
construct the weight vector $v$. One mthod is by regression.
Another is names the Fisher linear discriminant. A third is called
the
 support vector machine approach.

\subsubsection{Regression Classifier}

In the regression classifier, we form a matrix $A$ whose rows are
the feature vectors. We form a vector $b$ whose $kth$ component is
$0$ if the $kth$ feature vector comes from class one and whose
$kth$ component is $1$ is the $kth$ feature vector comes from
class two. We determine the weight vector $v$  as that vector that
minimizes $||Av-b||$. The minimizing vector $v$ is given by the
normal equation
$$
v=(A'A)^{-1}A'b
$$
The discriminant function is defined by
$$
f(x) = v'x
$$
We assign a vector to class one if $f(x) < \theta$ and to class
two otherwise. $\theta$ is chosen to minimize the error of the
assignment.

\subsubsection{Fisher Linear Discriminant}

Fisher's linear discriminant function is obtained by maximizing
the Fisher's discriminant ratio, which, as described below, is the
ratio of the projected between class scatter to the projected
within class scatter.

Let $v$ be the unknown weight vector. Let $\mu_1$ and $\mu_2$ be
the class one and two sample means and let $C_1$ and $C_2$ be the
class one and two sample covariance matrices. Let $N_1$ be the
number of feature vectors in class one and let $N_2$ be the number
of feature vectors in class two. Define the overall mean $\mu$ by
$$
\mu = P_1 \mu_1+P_2\mu_2
$$
where $P_1 = N_1/(N_1+N_2)$ and $P_2=N_2/(N_1+N_2)$. Then the
between-class scatter matrix $S_b$
 is given by
\begin{eqnarray*}
S_b &=& \sum_{i=1}^2 P_i ( \mu_i - \mu ) ( \mu_i - \mu )' \\
&=& P_1 P_2 (\mu_1 - \mu_2) ( \mu_1 - \mu_2)'
\end{eqnarray*}

Define $S_w$ to be the average class conditional scatter matrix,
then
$$S_w = \sum_{i=1}^2 P_i C_i$$
Finally, if we let $S$ designate the scatter matrix of the mixture
distribution,
$$S = \frac{1}{ N_1 + N_2} \sum_{k=1}^{N}
\left[ ( x_k - \mu ) ( x_k - \mu )' \right]$$ then
$$S = S_w + S_b$$

In the one dimensional projected space  one can easily show that
the projected between class scatter $s_b$ and the projected
within-class scatter $s_w$ are expressed as
$$s_b = v' S_b v$$
$$s_w = v' S_w v$$
Then the Fisher discriminant ratio is defined as

$$F(v) = \frac{s_b}{ s_w} = \frac{v' S_b v}{ v' S_w v}$$

The optimum direction $v$ can be found by taking the derivative of
$F(v)$ with respect to $v$ and setting it to zero:
$$\nabla F(v) = ( v' S_w v )^{-2} ( 2 S_b v v' S_w v - 2 v' S_b v S_w v ) = 0$$

\noindent From this equation it follows that
$$v' S_b v S_w v = S_b v v' S_w v$$
If we divide both sides by the quadratic term $v' S_b v$, then

\begin{eqnarray*}
S_w v &=& \frac {v' S_w v}{ v' S_b v} ) S_b v \\
&=& \lambda S_b v \\
&=& \lambda P_1 P_2 (\mu_1 - \mu_2) (\mu_1 - \mu_2)' v \\
&=& \lambda \kappa (\mu_1 - \mu_2)
\end{eqnarray*}

where $\lambda$ and $\kappa$ are some scalar values defined as
\begin{eqnarray*}
\lambda &=& \frac{v' S_w v}{ v' S_b v}\\
 \kappa &=&P_1 P_2 (\mu_1 - \mu_2)' v
\end{eqnarray*}

Thus we have the weighting vector $v$ as
$$v = K S_w^{-1} ( \mu_1 - \mu_2 )$$
where $K = \lambda \kappa$ is a multiplicative constant.

The discriminant function is defined by $f(x) = v'x$. The vector
$x$ is assigned to class one if $F(x) > \theta$ and to class two
otherwise. The threshold $\theta$ is set to a value that minimizes
the error of the class assignment.

\subsubsection{Support Vector}
Let $z_1,\ldots, z_N$ be the set of $N$ training vectors with
corresponding labels $y_1,\ldots, y_N$, a label being +1 for class
one and -1 for class two. Consider a hyperplane ribbon $R$ that
can separate the training vectors of class one from the training
vectors of class two. We represent $R$ by
$$
R=\{ x \ | \-1 \le w'x \le 1 \}
$$
The support vector machine approach seeks to find the widest
ribbon so that $R$ separates the vectors of class one from class
two.

The distance of the hyperplane $H=\{ x \ | \ w'x =1 \}$ from the
origin is $\frac{1}{||w||}$. Therefore the ribbon $R$ has width
$\frac{2}{||w|||}$. To maximize this width is equivalent to
minimizing $||w||$. This minimization must be done under the
constraint that $y_kw'z_k > 1$, k= 1,\ldots, K. Define the matrix
$A$ by
$$
A = \left(\begin{array}{c}
             y_1z'_1\\
             y_2z'_2\\
             \vdots\\
             y_Kz'_K\\
             \end{array}
             \right)
$$
Define the vector $b$ to be a vector of $K$ components all of
which have value $1$. The support vector approach determines the
vector $w$ by minimizing $w'w$ subject to the constraint $Aw > b$.
This can be solved by standard quadratic programming methods.

\subsection{Quantizing}
\label{se:quantizing}

 Let $f$ be a discriminant function. We
evaluate $f(x)$ over all the sampled vectors $x$ from class one
and from class two to determine the range. We divide the range in
a fixed number $M$ of quantizing intervals. The simplest way is
called equal interval quantizing. Here the range is divided up
into $M$ equal intervals. In each interval the number of sampled
vectors coming from class one and coming from class two is
determined. The interval is labelled by the class of the majority
of the vectors in it.

A vector $x$ having discriminant value $f(x)$ which falls into the
$mth$ quantizing interval is assigned to the class that labels the
quantizing interval.

Another simple alternative quantizing scheme is to divide the
range into intervals each of which have the same number of sampled
discriminant values. This is called equal probability quantizing.

A more complex scheme is to divide the discriminant range into $M$
intervals in such a way that the classification error is
minimized.

\subsection{Classification Trees}

The type of classification tree discussed  here is a binary tree
with a simple discriminant function; thus every nonterminal node
has exactly two children \cite{s:breiman}. During classification,
if the node's discriminant function is less than a threshold, the
left child is taken; if it is greater than the threshold, the
right child is taken. This section describes the design process of
the binary tree classifier using a simple discriminant function.
There are two methods of expanding a nonterminal node according to
the selection of a decision rule for the node. We show how to use
an {\it entropy purity function} to decide what the threshold
value should be, and we discuss the relationship of the purity
function to the $\chi^2$ test statistic. We discuss the criteria
for deciding when to stop expanding a node and for assigning a
class.

Let $x_1,\ldots,x_N$ be the set of vectors in the training set.
Associated with each vector is a class label.  Let $M$ be the
number of classes. Let \[X^n = \{ x^n_k \mid k = 1,\dots,N^n \}\]
be the subset of $N^n$ training vectors associated with node $n$.
Let $N^n_c$ be the number of training vectors for class $c$ in
node $n$. Since $N^n$ is the total number of training samples in
node $n$, we must have $N^n = \sum_{c=1}^{M} N^n_c$. The decision
rule selected for node $n$ is that discriminant function having
the greatest purity, a quality we will precisely define later.

Now we define how the decision rule works at node $n$. Consider
the feature vector $x^n_k$. If the discriminant function
$f(x^n_k)$ is less than or equal to the threshold, then $x^n_k$ is
assigned to class $\Omega^n_{LEFT}$, otherwise it is assigned to
class $\Omega^n_{RIGHT}$. An assignment to $\Omega^n_{LEFT}$ means
that the feature vector descends to the left child node. An
assignment to $\Omega^n_{RIGHT}$ means that the feature vector
descends to the right child node.

Given a discriminant function $f$, we sort the feature vectors in
the set $X^n$ in an ascending order according to their
discriminant function value. Without loss of generality we assume
that the feature vectors are sorted in such a way that $f(x^n_k)
\leq f(x^n_{k+1})$ for $k = 1, \dots, N^n-1,$. Let $w^n_k$ be the
true class associated with the measurement vector $x^n_k$. Then
the set of candidate thresholds $T^n$  is defined by
$$T^n = \left\{ = \frac{f(x^n_{k+1}) - f(x^n_k)} { 2} \Big| w^n_{k+1} \ne w^n_k
\right\}$$

For each possible threshold value, each feature vector
 $x^n_k$ is classified by using the decision rule
specified above. We count the number of samples $n^t_{Lc}$
assigned to $\Omega^n_{LEFT}$ whose true class is $c,$ and we
count the number of samples $n^t_{Rc}$ assigned to
$\Omega^n_{RIGHT}$ whose true class is $c$.
\begin{eqnarray*}
n^t_{L_c}&=&\#\{k\Big|f(x^n_k)\le t\hbox{ and }w^n_k=c\}\\
n^t_{R_c}&=&\#\{k\Big|f(x^n_k)>t\hbox{ and }w^n_k=c\}
\end{eqnarray*}

Let $n^t_L$ be the total number of samples assigned to
$\Omega^n_{LEFT}$ and $n^t_R$ be the total number of samples
assigned to $\Omega^n_{RIGHT}$, that is,
$$n^t_L = \sum_{c=1}^{M} n^t_{Lc}$$
$$n^t_R = \sum_{c=1}^{M} n^t_{Rc}$$

We define the purity $PR^t_n$ of such an assignment made by node
$n$ to be
$$PR_n = \sum_{c=1}^{M} \left( n^t_{Lc} \ln p^t_{Lc} +
n^t_{Rc} \log p^t_{Rc}\right)$$ where
$$p^t_{Lc} = \frac{n^t_{Lc}}{ n_L}$$
$$p^t_{Rc} = \frac{n^t_{Rc}}{ n_R}$$
The discriminant threshold selected is the threshold $t$ that
maximizes the purity value $PR^t_n$. The purity is such that it
gives maximum value when the classes of the training vectors are
completely separable. For example, consider a nonterminal node
having $m$ units in each of three classes in the training sample.
If the selected decision rule separates the training samples such
that the {\it LEFT} child contains all feature vectors
 in one class and the {\it RIGHT} child contains
all the feature vectors in the other two classes, the purity is $0
- 2 ( m \ln \frac{1}{ 2} ) = -2m\ln 2$. In the worst case, when
both the {\it LEFT} and the {\it RIGHT} children contain the same
number of feature vectors  for each class, the purity is $- 3 (
\frac{m}{ 2} \ln \frac{1}{ 3} ) - 3 ( \frac{m}{ 2} \ln \frac{1}{
3} ) = -3m \ln 3$. Thus we can easily see that the purity value of
the former case, where the training samples are completely
separable, is greater than the purity value of the latter case,
where the training samples are not separable.

The maximization of the purity can also be explained in terms of
the $\chi^2$ test of goodness of fit. If a decision node is
effective, the distribution of classes for the children nodes will
be significantly different from each other. A statistical test of
significance can be used to verify this. One test statistic that
measures the significance of the difference of the distributions
is defined by
$$\chi^2 = \sum_{c=1}^{M} \left(n^t_{Lc} \ln p^t_{Lc} +
 n^t_{Rc} \ln p^t_{Rc} -
 N^n_c \ln \frac{ N^n_c}{N^n }\right)$$
It has a $\chi^2$ distribution with $M-1$ degree of freedom.
Comparing this equation with $PR^t_n$, we find that the $\chi^2$
value is just the sum of the purity $PR^t_n$ and some constant
value.

Now we discuss the problem of when to stop the node expanding
process and how to assign a class to the terminal node. First, it
is not reasonable to generate a decision tree that has more
terminal nodes than the total number of training samples. Using
this consideration as a starting point, we set the maximum level
of the decision tree to be $\log_2 N^1 - 1$, which makes the
number of terminal nodes less than $\frac{N^1}{ 2}$, where $N^1$
is the number of training samples in node 1, the root node. Next,
if the $\chi^2$ value is small, the distributions of classes for
the children nodes are not significantly different from each
other, and the parent node need not be further divided. Finally,
when the number $N^n$ of units at node $n$ becomes small, the
$\chi^2$ test cannot give a reliable result. Therefore we stop
expanding node $n$ when $N^n$ is less than some lower limit. If
one of these conditions is detected, then the node $n$ becomes a
terminal node.

When a node becomes terminal, an assignment of a label to the node
is made. Each terminal node is assigned that class label that is
the majority of the class labels of the training vectors in the
node.

An alternative decision tree construction procedure uses the
probability of misclassification in place of entropy. In this
procedure the decision rule selected for the nonterminal node is
the one that yields the minimum probability of misclassification
of the resulting assignment. To describe the termination condition
of a node expansion, we first define type I and type II errors as
follows: Let type I error be the probability that a unit whose
true class is in $\Omega^n_{LEFT}$ is classified as
$\Omega^n_{RIGHT}$, and type II error be the probability that a
unit whose true class is in $\Omega^n_{RIGHT}$ is classified as
$\Omega^n_{LEFT}$. Then, if the sample space is completely
separable, we would get zero for both type I and type II errors.
Since this is not always the case, we control these errors by
considering only those thresholds in the process of threshold
selection that give type I error less than $\epsilon_I$ and type
II error less than $\epsilon_{II},$ where $\epsilon_I$ and
$\epsilon_{II}$ are values determined before we start constructing
the decision tree. Next, in the process of expanding a nonterminal
node, if we cannot find a decision rule that gives type I error
less than $\epsilon_I$ and type II error less than
$\epsilon_{II}$, which means that the sample space is not
separable at a $\epsilon_I$ and $\epsilon_{II}$ level, we stop
expanding this nonterminal node. This process of decision tree
construction is repeated until there is no nonterminal node left
or the level of the decision tree reaches the maximum level.
Assignment of a class to a nonterminal node is done in the same
way as in the previous procedure.

As just described, there are two groups of classes at each
nonterminal node in the binary decision tree. For the purpose of
discussion, we denote all the classes in one group as
$\Omega_{LEFT}$ and all the classes in the other group as
$\Omega_{RIGHT}$. The job of the decision rules discussed here is
to separate the left class $\Omega_{LEFT}$ from the right class
$\Omega_{RIGHT}$. We will employ the same notational conventions
used previously. The superscript $n$ denoting the node number will
be dropped from the expression if it is clear from the context
that we are dealing with one particular node $n$.

The simplest form for a discriminant function to take is a
comparison of one measurement component to a threshold.  This is
called a {\it threshold decision} rule. If the selected
measurement component is less than or equal to the threshold
value, then we assign class $\Omega_{LEFT}$ to the unit $u_k$;
otherwise we assign class $\Omega_{RIGHT}$ to it. This decision
rule requires a feature vector component index and a threshold.
Each feature vector component is selected in turn and the set of
threshold candidates $T$ of that component is computed. For each
threshold in the set $T$, all vectors in the training set $X^n$ at
node $n$ are classified into either class $\Omega_{LEFT}$ or class
$\Omega_{RIGHT}$ according to their value of the selected feature
component. The number of feature vectors for each class assigned
to class $\Omega_{LEFT}$ and to class $\Omega_{RIGHT}$ is counted,
and the entropy purity is computed from the resulting
classification. A threshold is selected from the set of threshold
candidates $T$ such that, when the set $X^n$ is classified with
that threshold, a maximum purity in the assignment results. This
process is repeated for all possible feature components, and the
component and threshold that yield an assignment with the maximum
purity is selected.

\subsection{Clustering}
\label{sec:kmeans}
All the previous methods required that the feature vectors be
generated for each of the objects to be discriminated and an
associated class label be associated with each of the feature
vectors. In this section we discuss a method to automatically
determine the natural classes, called clusters, to be associated
with each feature vector. Feature vectors with the same cluster
label are more similar to each other and less similar to feature
vectors of a different cluster label. Clustering can be used to
help generate hypotheses about natural distinctions between
mathematical objects before these distinctions themselves are
known.

The most widely used clustering scheme is called K-means. It is an
iterative method. It begins with a set of K cluster centers
$\mu^0_1,\ldots, \mu^0_K$. Initially these are chosen at random.
At iteration $t$ each feature vector is assigned to the cluster
center to which it is closest. This forms index sets
$S^t_1,\ldots, S^t_K$.

$$
S^t_k = \{ n \ | \ ||x_n - \mu^t_k|| \le ||x_n - \mu^t_m||,
m=1,\ldots,K \}
$$

Then each cluster center is redefined as the mean of the feature
vectors assigned to the cluster.

$$
\mu^{t+1}_k = \frac{1}{\#S^t_k} \sum_{n\in S^t_k} x_n
$$

Each iteration reduces the criterion function $J^t$.
$$
J^t = \sum_{k=1}^K \sum_{n\in S^t_k} ||x_n - \mu^t_k||
$$

As this criterion function is bounded below by zero, the
iterations must converge.

\section{Recognizing Whitehead Minimal Words in Free Groups}
In this section we return to the Whitehead minimal word problem.
We first discuss a pattern recognition system based on linear
regression,  then we show some applications of clustering.

\subsection{Linear regression classifier}
\label{subsec:generate}

\noindent
{\bf Data sets.}  To train a classifier, we must have a
training set. To test the classifier we must have an independent
test set.

A "random" element $w$  of  $F = F(X)$ can be produced as follows.
Choose   randomly   a number $l$ (the length of $w$), and a random
sequence $y_1, \ldots, y_l$ of elements $y_i \in X^{\pm 1}$ such
that $y_i \neq y_{i+1}^{-1}$, where $y_1$ is chosen randomly and
uniformly from $X^{\pm 1}$, and $y_{i+1}$ is  chosen randomly and
uniformly from the set $X^{\pm 1} - \{y_i^{-1}\}$. Similarly, one
can pseudo-randomly generate cyclically reduced words in $F$,
i.e., words $w = y_1 \ldots y_l$ where $y_1 \neq y_l^{-1}$.

 To generate the training data set  we used the following procedure.
 For each positive integer $l = 1, \ldots, 1000$ we generate  randomly
and  uniformly 10 cyclically reduced words from $F(X)$ of length
$l$. Denote the resulting set by  $W$. Then using the
deterministic Whitehead algorithm one can effectively construct
the corresponding set of minimal elements
\[W_{min} =\{w_{min} |  w \in W\}.\]
With probability 0.5 we substitute each $v \in W_{min}$ with the
word $\widetilde{v^t}$, where $t$ is a randomly and uniformly
chosen Whitehead automorphism such that $|\widetilde{v^t}|
> |v|$ (if $|\widetilde{v^t}| =  |v|$ we chose another automorphism $t$, and so
on). Now, the resulting set $D$ is a set of pseudo-randomly
generated cyclically reduced words representing the classes of
minimal and non-minimal elements  in approximately equal
proportions. We choose $D$ as the training set.

One remark is in order here. It seems, the class of non-minimal
elements in $D$ is not quite representative, since every one of its
elements $w$ has Whitehead complexity 1, i.e., there exists a
single Whitehead automorphism which reduces $w$ to $w_{min}$ (see
\cite{G:Miasnikov03GAWhitehead} for details on Whitehead
complexity). However, our experiments showed that the set $D$  is
a sufficiently good training dataset which is much easier to
generate than a set with uniformly distributed Whitehead
complexity of elements. A possible mathematical explanation of
this phenomena is mentioned in \cite{G:Miasnikov03GAWhitehead}.

To test and evaluate the pattern recognition methodology we
generate several test datasets of different type:
\begin{itemize}
 \item  A test set $S_e$ which is generated by the same procedure
as for the training set $D$, but independently of $D$.
 \item A test set $S_R$ of  randomly generated elements of $F(X)$.
   \item A test set $S_P$ of (pseudo-) randomly generated \emph{primitive}
elements in $F(X)$. Recall that $w \in F(X)$ is primitive if and
only if there exists
 a sequence of Whitehead automorphisms $t_1 \dots t_l \in \Omega(X)$
 such that $x^{ t_1 \dots
t_l} = w$ for some  $x \in X^{\pm 1}$.  Elements in $S_P$ are
generated by the  procedure described in
\cite{G:Miasnikov03GAWhitehead}, which, roughly speaking,  amounts
to a random choice of $x \in X^{\pm 1}$ and  a random choice of a
sequence of automorphisms $t_1 \dots t_l \in \Omega(X)$.
 \item A test set $S_{10}$ which is generated in a way
similar to the procedure used to generate the training set $D$.
The only difference is that the
 non-minimal elements are obtained by applying  not one, but several  randomly
chosen Whitehead automorphisms. The number of such automorphisms
is chosen uniformly randomly from the set $\{1, \ldots , 10\}$,
hence the name.
\end{itemize}

\bigskip
\noindent {\bf Features.} We use the feature vectors $f_1(w),
\ldots, f_6(w)$ described in Section \ref{se:feature}.

\bigskip
\noindent
 {\bf Model.} Our model is based on the linear
regression classifier described in Section
  \ref{se:linear}. For any word $w$ having feature vector $z(w)$ we compute the
discriminant function
$$
\hat P(w) = b'z(w)
$$
where $b'$ is the vector of regression coefficients obtained from
the training data set $D$. The decision rule is based on the equal
interval quantizing method described in Section
\ref{se:quantizing}.

\bigskip
\noindent
 {\bf Evaluation.} \label{se:conf}
 Let $D_{eval}$ be a test data set.   To evaluate the performance of the  given
{\bf PR} system we use a simple accuracy measure:
$$
A =|\{ w \in D_{eval} | \mbox{\rm minimality of}  \ w \ \mbox{\rm
is decided correctly} \}|
$$
It is the number of the correctly classified elements from the
test set $D_{eval}$.

\bigskip
\noindent
 {\bf Results.} \label{se:eval-results}
 Now  we present the
evaluation of  classifiers $\mathbf{P}_f$ on the test dataset
$S_e$ when $f$ runs over the set of feature mappings $f_1, \ldots,
f_6$ mentioned above. By $A(f)$ we denote the accuracy of the
classifier $\mathbf{P}_f$. For simplicity we present results only
for the free group $F$ of rank 2 with basis $X = \{a,b\}$.

  The results of evaluation of the classifiers $\mathbf{P}_i = \mathbf{P}_{f_i},
  i = 1, \ldots, 6$  on $S_e$ are given in Table \ref{tab:results}.

%
%
\begin{table}[ht]
\begin{center}
\begin{tabular}{|l|c|c|c|c|c|c|}
\hline
& $A(f_1)$ & $A(f_2)$ & $A(f_3)$ & $A(f_4)$ & $A(f_{5})$ & $A(f_{6})$ \\
\hline
$|w|>0$ &0.954&0.968&0.926&0.869&0.977&0.980 \\
\hline
$|w|>4$ &0.957&0.969&0.927&0.870&0.977&0.981 \\
\hline
$|w|>100$ &0.975&0.984&0.947&0.893&0.992&0.994 \\
\hline
\end{tabular}
\vskip 3ex
\end{center}
\caption{Performance of the classifiers  $\mathbf{P}_1, \ldots,
\mathbf{P}_6 $ on the set $S_e$.} \label{tab:results}
\end{table}
%

\begin{figure}[h]
\centerline{\includegraphics[scale=0.5]{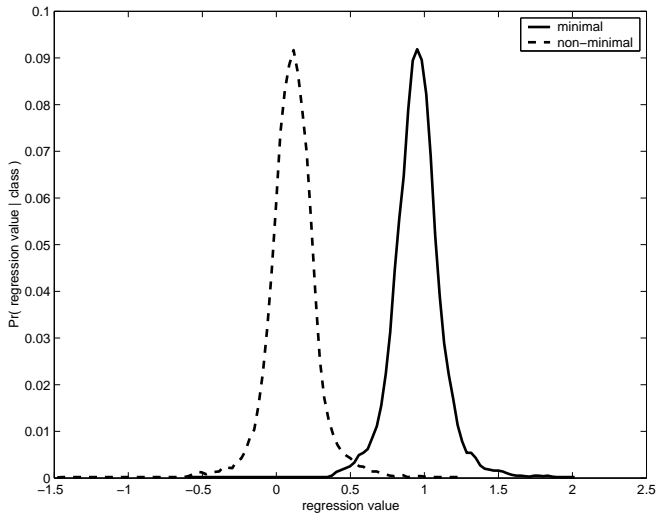}}
\caption{Conditional probabilities for $\mathbf{P}_{6}$.}
\label{fig:CondProbsf6}
\end{figure}

\bigskip
One can draw the following conclusions from the experiments:
\begin{itemize}
 \item  It seems, that the accuracy of the  classifiers increases when one adds
   new edges to the graphs related to the feature mappings (though it is not
   clear what  is the optimum set  of features);

\item  The classifier $\mathbf{P}_{6}$ is the best so far, it is
remarkably reliable;

 \item Very short words are difficult to classify (perhaps,
 because they do not provide sufficient information for the
 classifiers);

 \item  The estimated conditional probabilities for $\mathbf{P}_{6}$
 (which come from the Bayes' decision rule)  are presented
 in Figure \ref{fig:CondProbsf6}. Clearly,  the
classes of minimal and non-minimal elements are  separated around
0.5 with a small overlap. So the regression works perfectly with
the threshold $\Theta = 0.5$.
\end{itemize}

%
%

\bigskip
\noindent {\bf Looking for the best feature vectors.}
  As we have seen above the performance of a classifier $\mathbf{P}_{f}$
  directly depends on the
  feature vector $f(w)$  built into  it.   Sometimes it is possible to reduce the number of features in
$f$ maintaining the same level of classification accuracy of, and
even find more efficient combinations of the given features. The
corresponding procedure is called {\em feature selection}.

By far, the feature vector $f_6$ was the most effective. Observe,
that $f_6$ has 60 components (features). To find a better feature
vector we used an iterative greedy procedure  to select the best
vector from the  set of all counting functions of the type
\[\{ C(w,xvy) \mid x,y \in X^{\pm 1}, v \in F(X), 1 \leq |v| \leq 3 \}.\]
 It turns out that one of the most effective feature  vectors consists
only of two counting functions:
\[f^*(w) = \frac{1}{|w|}<C(w,a^{-1} b),C(w,b^{-1} a)>.\]

The
 results of comparison of $\mathbf{P}_{*} = P_{f^*}$ with $\mathbf{P}_{1}$ and $\mathbf{P}_{6}$
 are presented in Table \ref{tab:results-2}.
%
%
\begin{table}[h]
\begin{center}
\begin{tabular}{|l|c|c|c|c|}
\hline
& $A(f_1)$  & $A(f_{6})$ & $A(f^*)$  \\
\hline
$|w|>0$ &0.954&0.980&0.987 \\
\hline
$|w|>4$ &0.957&0.981&0.989 \\
\hline
$|w|>100$ &0.975&0.994&0.993 \\
\hline
\end{tabular}
\vskip 3ex
\end{center}
\caption{Comparative results for $\mathbf{P}_*$.}
\label{tab:results-2}
\end{table}
%

\subsection{Clustering}

In this section we describe one application of the K-means
clustering scheme to the Whitehead minimization problem.

In general cluster analysis is used to recover hidden structures
in a sampled set of objects. In the  Whitehead's  minimization
 method the following {\em Length Reduction Problem}  is of prime interest:
 given a non-minimal
 word
 $w \in F$ find a ({\em length-reducing})
 Whitehead automorphism $t$ such that
 \[|wt| < |w|.\]

Below we apply K-means clustering method described in section
\ref{sec:kmeans} to the Length Reduction Problem. The task is to
partition a given set of non-minimal elements into clusters in
such a way that every cluster would have  a Whitehead automorphism
assigned to it which reduces the length of the most words from the
cluster. To illustrate performance of the K-means algorithm we
take the feature vector function  $f_2$ from Section
\ref{se:feature} and the standard Euclidean  metric in
$\mathbb{R}^4$ (recall that $f_2(w) \in \mathbb{R}^4$ for $w \in
F_2$).

To perform the K-means procedure one needs to specify in advance
the number K of expected clusters. Since we hope that every such
cluster $\mathcal{C}$ will correspond to a particular Whitehead
automorphism that reduces the length of elements in $\mathcal{C}$
then the expected number of clusters can be calculated as follows.
 It is easy to see that  the set $\Omega_2$ of all Whitehead automorphisms
 of the free group  $F_2(X)$ with basis  $X = \{a,b\}$ splits
 into two subsets: the set
  $$N_2 = \left \{ \TwoMap{ab}{b}, \TwoMap{b^{-1}a}{b}, \TwoMap{a}{ba},
  \TwoMap{a}{a^{-1}b} \right \}$$
of  Nielsen automorphisms and the set of conjugations. If we view
elements of $F_2$ as cyclic words (i.e., up to a cyclic
permutation) then the conjugations from $\Omega_2$ can be ignored
in the length reduction problem. Therefore, we would like the
K-means algorithm to find precisely 4 clusters.

 Let $S \subset F_2$
be a set of  non-minimal cyclically reduced words from $F_2$. We
construct the set
$$D = <f_2(w) \mid w \in S>$$
of feature vectors, corresponding to words in $S$.
 To start the K-means algorithm one needs to choose the set of initial centers
$\mu^0_t, t \in N_2$. Observe, that the algorithm is quite
sensitive to this choice of the centers. There are various methods
to generate $\mu^0_t, t \in N_2$, here we describe just one of
them. Let $S^\prime$ be a sample subset of the set $S$.
  For an automorphism $t \in N_2$ put
 $$C_t = \{ w \in S^\prime \mid |wt| < |w|, \forall r \in N_2(r \neq t \rightarrow
  |wr| \geq |w|) \}$$
 and define
\begin{equation}
\label{eq:centers}
\mu^0_t = \frac{1}{|C_t|} \sum_{w \in C_t} f_2(w)
\end{equation}
as the initial estimates for the cluster centers (we assume here
that the sets $C_t$ are not empty).

The goodness of the clustering is evaluated using a measure
$R_{max}$ defined below.
 Let  $\mathcal{C} \subset D$ be a cluster. For $t \in N_2$ define
$$R(t,\mathcal{C}) = \frac{|<v | |t(v)|<|v|, v \in \mathcal{C}>|}{|\mathcal{C}|}.$$
 The number
$R(t,\mathcal{C})$ shows how many elements in $\mathcal{C}$ are
reducible by $t$. Now put
$$R_{max}(\mathcal{C}) = \max \{R(t,\mathcal{C}) \mid  t \in N_2 \}.$$
The number $R_{max}(\mathcal{C})$ shows how many elements in
$\mathcal{C}_i$ can be reduced by a single automorphism.

Now we perform 4-means clustering algorithm on the set $S$  and
compute $R(t,\mathcal{C}_i)$, $t \in N_2$,  for each obtained
cluster $\mathcal{C}_i$. The results of the experiments are given
in Table \ref{tab:kmeans} for two different choices of the initial
centers (the  random choice and the choice according to
\ref{eq:centers}). We can see from the table that the clustering,
indeed, groups elements in $S$ with respect to the length reducing
transformation.

\begin{table}[ht]
\begin{center}
\begin{tabular}{|l|c|c|}
\hline
 & random $\mu^0_t$& $\mu^0_t$ estimated by (\ref{eq:centers}) \\
\hline
avg$(R_{max})$ & 0.930 & 1.000 \\
\hline
max$(R_{max})$ &1.000 & 1.000 \\
\hline
min$(R_{max})$ & 0.743& 1.000 \\
\hline
\end{tabular}
\vskip 3ex \label{tab:kmeans} \caption{ Evaluation of 4-means
clustering of the set $D$. }
\end{center}
\end{table}

The experiments show that it is possible to cluster non-minimal
elements in $F_2$, using the standard clustering algorithms,
 in such a way  that:
  \begin{itemize}
  \item  every cluster
contains elements whose length  can be (with a very high
probability) reduced by a particular Nielsen transformation;
 \item the  transformation that reduces the length of
the most elements from one cluster does not reduce the length of
the most elements in another cluster.
 \end{itemize}
 This gives a very strong heuristic for choosing a length reducing
automorphism for a given word $w \in F_2$.

A simple decision rule which for a given word $w \in F_2$ will
predict a corresponding length reducing automorphism $t$ can be
defined as follows. Let $\mu_t$ be the centers  of the clusters,
produced by the K-means clustering. Each $\mu_t$ corresponds to a
particular automorphism $t \in N_2$. For a given non-minimal
cyclically reduced word $w \in F_2$ we select an automorphism $t^*
\in N_2$ such that
\[  \forall t \in N_2 \ \left (||f_2(w) - \mu_{t^*}|| \ \leq \ ||f_2(w) - \mu_t|| \right ) \]
 as the most probable length-reducing automorphism for $w$.

In the conclusion we would like to add that a similar analysis can
be used to predict most probable length-reducing automorphisms for
words in free groups of ranks $n$ larger then 2. However, the
number of the corresponding clusters grows exponentially with $n$
which increases the error rate of the classification. In this case
more careful clustering still could be applied where  the clusters
correspond to some particular groups of Whitehead automorphism.

\bibliographystyle{abbrv}
\bibliography{ThesisBibliography}
\end{document}